\documentclass[12pt]{amsart}
\theoremstyle{plain}
\newtheorem{Thm}{Theorem}

\newtheorem{Cor}[Thm]{Corollary}

\newtheorem{Rk}[Thm]{Remark}

\begin{document}

\title[a non-local population model]
{A non-local population model of logistic type equation}

\author{Li Ma, Liang Cheng}

\address{Department of mathematical sciences \\
Tsinghua university \\
Beijing 100084 \\
China} \email{lma@math.tsinghua.edu.cn} \dedicatory{}
\date{September 26th, 2009}

\begin{abstract}

In this paper, we propose a new non-local population model of
logistic type equation on a bounded Lipschitz domain in the whole
Euclidean space. This model preserves the $L^2$ norm, which is
called mass, of the solution on the domain. We show that this
model has the global existence, stability and asymptotic behavior
at time infinity.


\end{abstract}

\thanks{$^*$ The research is partially supported by the National Natural Science
Foundation of China 10631020 and SRFDP 20060003002. }
 \maketitle

 \section{Introduction}
 In this work, we propose a new population model with non-local
 term. This non-local term keeps the mass, a global quantity to be defined below,
 constant for all time. Let's first review some previous study of
 population modelings.

After a critical study of Malthus's population model, people
believe that a good population model should have good behavior.
This makes it come to the logistic model, which is a slight
modification of Malthus's model. By definition, the logistic model
is a population model such that it describes the changes over time
of a population occupying a single small region. For more models
and the history of population modeling, we refer to the work
\cite{U04}, and the books \cite{K95} and \cite{R91}. In
mathematical language, the logistic model can be stated as below.
Let $P$ be the population quantity. Then the change rate of $P$ is
the difference between the birth rate $\frac{d B}{dt}$ and the
death rate $\frac{d D}{dt}$, i.e.,
$$
\frac{d P}{dt}=\frac{d B}{dt}-\frac{d D}{dt}.
$$
From the experimental observation, we put
$$
\frac{d B}{dt}=aP+bP^2
$$
and
$$
\frac{d D}{dt}=cP+dP^2.
$$
where $a,b,c,d $ are constants such that $a>c$ and $d>b$. Hence we
obtain
$$
\frac{d P}{dt}=(a-c)P-(d-b)P^2.
$$
Let
$$
r=a-c
$$
and
$$
K=\frac{a-c}{d-b}.
$$
Then we get the logistic model
$$
\frac{d P}{dt}=rP(1-\frac{P}{K}).
$$
 Here the  growth rate $r$ represents the change at which the population may grow
 if it were unencumbered by environmental degradation, and
the parameter $K$ represents the carrying capacity of the system
considered. By definition, the carrying capacity is the population
level at which the birth and death rates of a species exactly
match, resulting in a stable population over time. Hence in some
modeling, one may assume that $K=K(t)$ is a periodic function in
time variable.

 The drawback of the model above is
that it ignores the impact of the environmental condition to the
population. When the environmental condition on the region $D$, a
bounded Lipschitz domain in $R^n$,is considered, one encounters
the following diffusion model of logistic type on $D$
\begin{equation}\label{eq1}
u_t=\Delta u+ru(1-\frac{u}{K}),
\end{equation}
where $u=u(t,x)$ is the population quantity such that $u=u(t,x)>0$
for $x\in D$ and $u(t,x)=0$ on $\partial D$, $t>0$. (\ref{eq1})
will be also called the Logistic equation as considered in
\cite{DM}. Note that the equation (\ref{eq1}) is a local model
such that the value $u(x,t)$ at $(x,t)$ depends only on its
immediate surroundings. We refer to \cite{DM}, \cite{Jiang}, and
\cite{shi09} for related models.

We now turn to our main subject. We now turn to our new subject.
In this paper, we shall modify (\ref{eq2}) to the following
non-local logistic equation
\begin{equation*}
\left\{
\begin{array}{ll}
         \partial_t u=\Delta u+\lambda(t)u+a(x)(u-u^p) \ \ \ \text{in}\ D\times\mathbb{R}_{+}, \\
          u(x,0)=g(x) \ \ \ \text{in}\ D,\\
          u(x,t)=0, \ \ \ \text{on}\ \partial D
\end{array}
\right.
\end{equation*}
where $p>1$ and $a(x)>0$ is a non-trivial Lipschtiz function on
the closure of the domain $D$, which has the positive solution and
preserves the $L^2$ the norm. By definition, we call the integral
quantity of $u(x,t)$, $\int_D u^2dx$, the mass of the population
model. Likewise,
$$\frac{1}{2}\frac{d}{dt}\int_{D}u^2dx=\int_{D}u
u_t=-\int_{D}|\nabla
u|^2dx+\lambda(t)\int_{D}u^2dx+\int_{D}a(u^2-u^{p+1})dx.$$ Thus,
one must have $\lambda(t)=\frac{\int_{D}(|\nabla
u|^2+a(u^{p+1}-u^2))dx}{\int_{D}g^2dx}$ to preserve the $L^2$
norm. Without loss of generality we assume $\int_{D}g^2dx=1$. Then
we consider the following problem on the bounded domain $D$
\begin{equation} \label{eq2}
\left\{
\begin{array}{ll}
         \partial_t u=\Delta u+\lambda(t)u+a(u-u^p)\ \ \ \text{in}\ D\times\mathbb{R}_{+} \\
          u(x,0)=g(x) \ \ \ \text{in}\ D \\
          u(x,t)=0 \ \ \ \text{on}\ \partial D
\end{array}
\right.
\end{equation}
where $p>1$, $\lambda(t)=\int_{D}(|\nabla u|^2+a(u^{p+1}-u^2)) dx$,
$g(x)\geq 0\ \text{in}\ D$, $\int_{D}g^2dx=1$ and $g\in C^1(D)$.

Similar to the global existence results obtained in C.Caffarelli
and F.Lin \cite{CL09} and our previous work \cite{MC2} (see also
related works \cite{LM}, \cite{MC}, and \cite{MA}), we have
following result.

\begin{Thm}\label{thm2}
Problem (\ref{eq2}) has a global solution $u(t)\in
L^{\infty}(\mathbb{R}_{+},H^1_0(D))\cap
L^{\infty}(\mathbb{R}_{+},L^{p+1}(D))\cap
L_{loc}^2(\mathbb{R}_{+},H^2(D))$.
\end{Thm}

\begin{Rk}
We note that solutions of (\ref{eq2}) have automatically higher
regularity for $t> 0$. Indeed, the bound of $\lambda(t)$ (see
(\ref{estimate11})) and the standard parabolic estimates imply that
solutions are H\"{o}lder continuous. Then coming back to
$\lambda(t)$, it would be a H\"{o}lder continuous function in time.
A bootstrap argument implies that $u$ is smooth in both spatial and
time variables if we assume $a$ is smooth function.
\end{Rk}

We also have the stability results for (\ref{eq2}).

\begin{Thm}\label{thm6}
Let $u,v$ be the two bounded solutions to problem (\ref{eq2}) with
initial data $g_u,g_v$ at $t=0$, where $g_u,g_v\in H^1(D)\cap
L^{\infty}(D)$. Then
$$
||u-v||_{L^2}^2\leq ||g_u-g_v||_{L^2}^2
\exp(C_1 t)
$$
and
$$
||u-v||_{H^1}^2\leq ||g_u-g_v||_{H^1}^2 \exp(C_2 t),
$$ where $C_1,C_2$ are the constants depending on the upper bound
of $||g_u||_{H^1(D)}, ||g_v||_{H^1(D)}$ and $||g_u||_{L^{\infty}},
||g_v||_{L^{\infty}}$. In particular, the solution to problem
(\ref{eq2}) is unique.
\end{Thm}

As the simple applications to theorem \ref{thm2}, we can study
asymptotic behavior of $u(t)$ of problem (\ref{eq2}).

\begin{Cor}\label{cor2}
Suppose $u(t)$ is the solution to problem (\ref{eq2}). Then one
can take $t_i\to \infty$ such that
$\lambda(t_i)\to\lambda_{\infty}$, $u(x,t_i)\rightharpoonup
u_{\infty}(x)$ in $H^1_0(D)$ and $u_{\infty}$ solves the equation
$\Delta
u_{\infty}+\lambda_{\infty}u_{\infty}+a(u_{\infty}-u_{\infty}^p)=0$
in $D$ with $\int_{D}|u_{\infty}|^2dx=1$.
\end{Cor}

In the next section, we prove
 the global existence, stability and asymptotic behavior of
solutions to the problem (\ref{eq2}). In particular, we give the
proofs of theorem \ref{thm2} to theorem \ref{thm6}, and corollary
\ref{cor2}. In section \ref{sect3} we give our conclusion based on
our study of problem (\ref{eq2}).

\section{global existence and stability  property}\label{sect2}
In this section we study the global existence, stability and
asymptotic behavior of solutions to the problem (\ref{eq2}).

\textbf{Proof of theorem \ref{thm2}.} Let us define a series
$u^{(k)}$ as
\begin{eqnarray} \label{eq4}
\left\{
\begin{array}{ll}
         u^{(0)}=g, \lambda^{(k)}(t)=\int_{D}(|\nabla u^{(k)}|^2+a(x)((u^{(k)})^{p+1}-(u^{(k)})^{2}) dx, \\
          \partial_t u^{(k+1)}=\Delta u^{(k+1)}+ \lambda^{(k)}(t)u^{(k+1)}+a(x)(u^{(k+1)}-(u^{(k+1)})^{p}),        \\
          u^{(k+1)}(x,0)=g(x),\\
          u^{k+1}(x,t)=0, \ \ \ \text{on}\ \partial D
\end{array}
\right.
\end{eqnarray}
a series of initial boundary value problems of linear parabolic
systems.

To prove the convergence of series $\{u^{(k)}\}$ constructed above,
we estimate for $k\geq 0$
\begin{equation}\label{bbb1}
\frac{1}{2}\frac{d}{dt}\int_{D}|\nabla u^{(k+1)}|^2
dx+\int_{D}|\Delta u^{(k+1)}|^2 dx+p\int_{D}
a(x)(u^{(k+1)})^{p-1}|\nabla u^{(k+1)}|^2dx
\end{equation}
$$
=\lambda^{(k)}(t)\int_{D}|\nabla u^{(k+1)}|^2 dx
+\int_{D}a(x)|\nabla u^{(k+1)}|^2 dx+\int_{D}(\nabla u^{(k+1)}\cdot
\nabla a)(u^{(k+1)}-(u^{(k+1)})^p)dx,
$$
\begin{equation}\label{ccc1}
\frac{1}{2}\frac{d}{dt}\int_{D}|\nabla u^{(k+1)}|^2
dx+\int_{D}|u_t^{(k+1)}|^2
dx+\frac{1}{p+1}\frac{d}{dt}\int_{D}a(x)(u^{(k+1)})^{p+1}dx
\end{equation}
$$
=\frac{\lambda^{(k)}(t)}{2}\frac{d}{dt}\int_{D}|u^{(k+1)}|^2 dx
+\frac{1}{2}\frac{d}{dt}\int_{D}a(x)|u^{(k+1)}|^2 dx.
$$
\begin{equation}\label{ddd1}
\frac{1}{p+1}\frac{d}{dt}\int_{D}(u^{(k+1)})^{p+1}
dx+\int_{D}p(u^{(k+1)})^{p-1}|\nabla
u^{(k+1)}|^2dx+\int_{D}a(x)(u^{(k+1)})^{2p}dx
\end{equation}
$$
=\lambda^{(k)}(t)\int_{D}(u^{(k+1)})^{p+1}dx
+\int_{D}a(x)(u^{(k+1)})^{p+1}dx.
$$
Now we denote $M=||a||_{C^1(D)}$. By (\ref{bbb1}), we get
\begin{eqnarray}\label{ineq11}
&&\frac{1}{2}\frac{d}{dt}\int_{D}|\nabla u^{(k+1)}|^2 dx\\
&\leq&(\lambda^{(k)}(t)+M)\int_{D}|\nabla u^{(k+1)}|^2
dx+\int_{D}(\nabla u^{(k+1)}\cdot
\nabla a)(u^{(k+1)}-(u^{(k+1)})^p)dx\nonumber\\
&\leq&(\lambda^{(k)}(t)+M)\int_{D}|\nabla u^{(k+1)}|^2 dx
+\frac{1}{4\epsilon}\int_{D}|\nabla u^{(k+1)}|^2|\nabla
a|^2dx\nonumber\\
&&+\epsilon\int_{D}(u^{(k+1)}-(u^{(k+1)})^p)^2dx\nonumber\\
&\leq&(\lambda^{(k)}(t)+M +\frac{M^2}{4\epsilon})\int_{D}|\nabla
u^{(k+1)}|^2dx+\epsilon c_1\int_{D}(u^{(k+1)})^{2p}dx,\nonumber
\end{eqnarray}
where $c_1$ is a constant only depending on $D$. Moreover, by
(\ref{ddd1}), we have
\begin{equation}\label{ineq12}
\frac{1}{p+1}\frac{d}{dt}\int_{D}(u^{(k+1)})^{p+1}
dx+\int_{D}a(x)(u^{(k+1)})^{2p}dx
\leq(\lambda^{(k)}(t)+M)\int_{D}(u^{(k+1)})^{p+1}dx.
\end{equation}
Now, we denote $\widetilde{\lambda}(t)=\int_{D}|\nabla
u|^2dx+M\int_{D}u^{p+1}dx$. Note that $a(x)$ is a positive
Lipschitz function on the compact domain $D$ and we may assume
$a(x)\geq c_0>0$. Hence
 by (\ref{ineq11}) and (\ref{ineq12}), we have
 \begin{eqnarray}\label{ineq13}
\frac{d}{dt}\widetilde{\lambda}^{k+1}(t)+((p+1)Mc_0-2\epsilon
c_1)\int_{D}(u^{(k+1)})^{2p}dx\leq
c_2(\lambda^{k}(t)+c_2)\widetilde{\lambda}^{k+1}(t),\nonumber
\end{eqnarray}
where $c_2$ is a constant only depending on $p,M,\epsilon$. We
choose $\epsilon$ such that $(p+1)Mc_0=3\epsilon c_1$ and denote
$c_3=\epsilon c_1$, combining with the fact $\lambda^{k}(t)\leq
\widetilde{\lambda}^{k}(t)$, we get
 \begin{eqnarray}\label{ineq14}
\frac{d}{dt}\widetilde{\lambda}^{k+1}(t)+c_3\int_{D}(u^{(k+1)})^{2p}dx\leq
c_2(\widetilde{\lambda}^{k}(t)+c_2)\widetilde{\lambda}^{k+1}(t).
\end{eqnarray}
Hence
 \begin{eqnarray}\label{ineq15}
  \\
\widetilde{\lambda}^{k+1}(t)\leq (\int_{D}|\nabla
g|^2dx+M\int_{D}g^{p+1}dx+c_2)\exp(c_2\int^t_0
(\widetilde{\lambda}^{k}+c_2)dt)\nonumber.
\end{eqnarray}
By induction, there is $\delta$ depending only on $\int_{D}|\nabla
g|^2dx$, $M\int_{D}g^{p+1}dx$ and $c_2$ such that
\begin{equation}\label{estimate11}
\lambda^{(k+1)}(t)\leq \widetilde{\lambda}^{k+1}(t)\leq c_4,\ \
\text{for}\ t\in[0,\delta], k\geq1,
\end{equation}
where $c_4$ is a constant depending on $\int_{D}|\nabla g|^2dx$,
$M\int_{D}g^{p+1}dx$ and $c_2$. Hence
\begin{equation}\label{estimate12}
\int_{D}|\nabla u^{(k+1)}|^2dx\leq c_4,\ M\int_{D}|
u^{(k+1)}|^{p+1}dx\leq c_4 \ \ \text{for}\ t\in[0,\delta], k\geq1.
\end{equation}
Integrate (\ref{ineq14}) with t, we can conclude that
\begin{equation}\label{ineq16}
\int^{\delta}_0\int_{D}|u^{(k+1)}|^{2p} dxdt\leq c_5.
\end{equation}
Now integrate (\ref{bbb1}) with t, also by (\ref{ineq16}), we get
$$
\frac{1}{2}\int_{D}|\nabla u(\delta)^{(k+1)}|^2
dx-\frac{1}{2}\int_{D}|\nabla u(0)^{(k+1)}|^2
dx+\int^{\delta}_0\int_{D}|\Delta u^{(k+1)}|^2 dxdt
$$
$$
+\int^{\delta}_0\int_{D}p(u^{(k+1)})^{p-1}|\nabla
u^{(k+1)}|^2dxdt\leq\int^{\delta}_0(\lambda^{(k)}(t)+M+\frac{M^2}{4\epsilon})\int_{D}|\nabla
u^{(k+1)}|^2 dxdt+c_3c_5.
$$
Hence
\begin{equation}\label{estimate13}
\int^{\delta}_0\int_{D}|\Delta u^{(k+1)}|^2 dxdt\leq c_6,
\end{equation}
where $c_6$ depending on $\int_{D}|g|^2dx$, $\int_{D}|\nabla g|^2
dx$, $c_3,c_5$, $M$ and $\delta$. Integrate (\ref{ccc1}) with t, we
get
\begin{eqnarray*}
&& \frac{1}{2}\int_{D}|\nabla u(\delta)^{(k+1)}|^2
dx-\frac{1}{2}\int_{D}|\nabla u(0)^{(k+1)}|^2
dx+\int^{\delta}_0\int_{D}|u_t^{(k+1)}|^2 dxdt\\
&&+
\frac{1}{p+1}\int_{D}a(x)(u^{(k+1)}(t))^{p+1}dx-\frac{1}{p+1}\int_{D}a(x)g^{p+1}dx\\
&=&
\int^{\delta}_0\frac{\lambda^{(k)}(t)}{2}\frac{d}{dt}\int_{D}|u^{(k+1)}|^2
dxdt +\frac{1}{2}\int_{D}a(x)|u(\delta)^{(k+1)}|^2
dx-\frac{1}{2}\int_{D}a(x)|u(0)^{(k+1)}|^2 dx.
\end{eqnarray*}
Hence
\begin{equation}\label{estimate14}
\int^{\delta}_0\int_{D}|u_t^{(k+1)}|^2 dxdt\leq c_7,
\end{equation}
where $c_7$ depending on $\int_{D}|g|^2dx$, $\int_{D}|\nabla g|^2dx$
, $\int_{D}g^{p+1}dx$, $c_3,c_5$, $M$ and $\delta$.

By (\ref{estimate11}), (\ref{estimate12}), (\ref{estimate13}) and
(\ref{estimate14}), there is a subsequence of $\{u^{(k)}\}$ (still
denoted by $\{u^{(k)}\}$) and a function $u(t)\in
L^{\infty}([0,\delta],H^1(D))\cap L^2([0,\delta],H^2(D))\cap
L^{\infty}([0,\delta],L^{p+1}(D))$ with $\partial_t u(t)\in
L^2([0,\delta],L^2(D))$ such that $u^{(k)}\rightharpoonup u$
weak$^\ast$ in $L^{\infty}([0,\delta],H^1(D))$, weakly in
$L^2([0,\delta],H^2(D))$ and weakly in
$L^{\infty}([0,\delta],L^{p+1}(D))$. Then we have $u^{(k)}\to u$
strongly in $L^2([0,\delta],H^1(D))$ and $u(t)\in
C([0,\delta],L^2(D))$. Hence $\lambda^{(k)}(t)\to \lambda(t)$
strongly in $L^2([0,\delta])$. Thus, we get a local strong solution
to problem (\ref{eq2}). Next, starting from $t=\delta$ we can extend
the local solution to $[0,2\delta]$ in exactly the same way as
above. By induction, we have a global solution to problem
(\ref{eq2}).
 $\Box$

The stability result will be proved in the similar manner as in
\cite{MC2}.

 \textbf{Proof of theorem \ref{thm6}.}
 By the arguments in theorem
\ref{thm2}, we can take a constant $C$ such that all
$||u||_{L^{\infty}(\mathbb{R}_{+},H^1(D))}$,
$||v||_{L^{\infty}(\mathbb{R}_{+},H^1(D))}$,
$||u||_{L^{\infty}(\mathbb{R}_{+},L^{\infty}(D))}$,
$||v||_{L^{\infty}(\mathbb{R}_{+},L^{\infty}(D))}$,
$||\lambda_u(t)||_{L^{\infty}(\mathbb{R}_{+})}$ and
$||\lambda_v(t)||_{L^{\infty}(\mathbb{R}_{+})}$ not less than $C$,
where $C$ is only depending on upper bound of $||g_u||_{H^1(D)},
||g_v||_{H^1(D)}$, $||g_u||_{L^{\infty}(D)},
||g_v||_{L^{\infty}(D)}$. We still denote $M=||a||_{W^{1,\infty}}$.
First we calculate
\begin{eqnarray*}
\frac{1}{2}\frac{d}{dt}\int_{D}(u-v)^2dx
&=&\int_{D}(u-v)(u_t-v_t)dx\\
&=&\int_{D}(u-v)(\Delta(u-v)+\lambda_u(t)u-\lambda_v(t)v+a(x)(u-v)\\
&&-a(x)(u^p-v^p))dx\\
&\leq&-\int_{D}|\nabla(u-v)|^2dx+\int_{D}(u-v)(\lambda_u(t)u-\lambda_v(t)v)dx\\
&&+\int_{D}a(x)(u-v)^2dx
\end{eqnarray*}
Note that
\begin{eqnarray*}
&&\int_{D}(u-v)(\lambda_u(t)u-\lambda_v(t)v)dx\\
&=&(\lambda_u(t)-\lambda_v(t))\int_{D}(u-v)udx+\lambda_v(t)\int_{D}(u-v)^2dx\\
&\leq&|\lambda_u(t)-\lambda_v(t)|(\int_{D}(u-v)^2dx)^{\frac{1}{2}}(\int_{D}u^2dx)^{\frac{1}{2}}
+\lambda_v(t)\int_{D}(u-v)^2dx\\
&\leq&C|\lambda_u(t)-\lambda_v(t)|(\int_{D}(u-v)^2dx)^{\frac{1}{2}}
+C\int_{D}(u-v)^2dx,
\end{eqnarray*}
and
\begin{eqnarray}\label{lambda1}
&&|\lambda_u(t)-\lambda_v(t)|\\
&=&|\int_{D}((|\nabla u|^2-|\nabla v|^2)+a(u^{p+1}-v^{p+1}))-a(u^2-v^2)dx|\nonumber\\
&\leq&\int_{D}|\nabla(u-v)|(|\nabla u|+|\nabla
v|)dx+|\int_{D}a(u-v)(\frac{u^{p+1}-v^{p+1}}{u-v})dx|\nonumber\\
&&+\int_{D}a(u-v)(u+v)dx\nonumber\\
&\leq&C(\int_{D}|\nabla(u-v)|^2dx)^{\frac{1}{2}}+MC(\int_{D}(u-v)^2dx)^{\frac{1}{2}}\nonumber.
\end{eqnarray}
We have
\begin{eqnarray*}
&&\frac{1}{2}\frac{d}{dt}\int_{D}(u-v)^2dx\\
&\leq&-\int_{D}|\nabla(u-v)|^2dx+C^2(\int_{D}|\nabla(u-v)|^2dx)^{\frac{1}{2}}(\int_{D}(u-v)^2dx)^{\frac{1}{2}}\\
&&+(MC^2+M+C)\int_{D}(u-v)^2dx\\
&\leq&-\frac{1}{2}\int_{D}|\nabla(u-v)|^2dx+(\frac{C^4}{2}+MC^2+M+C)\int_{D}(u-v)^2dx.
\end{eqnarray*}
By the Gronwall inequality \cite{Evans}, we have
$$
||u-v||_{L^2}^2\leq ||g_u-g_v||_{L^2}^2
\exp((\frac{C^4}{2}+MC^2+M+C)t).
$$
Further more,
\begin{eqnarray*}
&&\frac{1}{2}\frac{d}{dt}\int_{D}|\nabla(u-v)|^2dx\\
&=&-\int_{D}\Delta(u-v) \cdot (u-v)_tdx\\
&=&-\int_{D}\Delta(u-v)\cdot(\Delta(u-v)+\lambda_u(t)u+a(u-v)-a(u^p-v^p))dx\\
&=&-\int_{D}(\Delta(u-v))^2dx+\int_{D}\nabla(u-v)\cdot\nabla(\lambda_u(t)u-\lambda_v(t)v)dx\\
&&+\int_{D}\nabla(u-v)\cdot\nabla(a(u-v))dx+\int_{D}\Delta(u-v)\cdot
a(u^p-v^p)dx.
\end{eqnarray*}
Note that
\begin{eqnarray*}
&&\int_{D}\Delta(u-v)\cdot a(u^p-v^p)dx\\
&\leq&\frac{1}{2}\int_{D}(\Delta(u-v))^2dx+\frac{1}{2}\int_{D}a^2(u^p-v^p)^2dx\\
&=&\frac{1}{2}\int_{D}(\Delta(u-v))^2dx+\frac{1}{2}\int_{D}a^2(u-v)^2(\frac{u^p-v^p}{u-v})^2dx\\
&\leq&\frac{1}{2}\int_{D}(\Delta(u-v))^2dx+\frac{M^2C}{2}\int_{D}(u-v)^2dx
\end{eqnarray*}
Likewise,
\begin{eqnarray*}
&&\int_{D}\nabla(u-v)\cdot\nabla(\lambda_u(t)u-\lambda_v(t)v)dx\\
&=&(\lambda_u(t)-\lambda_v(t))\int_{D}\nabla(u-v)\cdot \nabla udx+\lambda_v(t)\int_{D}|\nabla(u-v)|^2dx\\
&\leq&|\lambda_u(t)-\lambda_v(t)|(\int_{D}|\nabla(u-v)|^2dx)^{\frac{1}{2}}(\int_{D}|\nabla
u|^2dx)^{\frac{1}{2}}
+\lambda_v(t)\int_{D}|\nabla(u-v)|^2dx\\
&\leq&C|\lambda_u(t)-\lambda_v(t)|(\int_{D}|\nabla(u-v)|^2dx)^{\frac{1}{2}}
+C\int_{D}|\nabla(u-v)|^2dx\\
&\leq&(C^2+\frac{MC^2}{2}+C)\int_{D}|\nabla(u-v)|^2dx+\frac{MC^2}{2}\int_{D}(u-v)^2dx,
\end{eqnarray*}
where the second inequality follows by (\ref{lambda1}). Furthermore,
we calculate
\begin{eqnarray*}
&&\int_{D}\nabla(u-v)\cdot\nabla(a(u-v))dx\\
&=&\int_{D}a|\nabla(u-v)|^2dx+\int_{D}(u-v)(\nabla(u-v)\cdot\nabla
a)dx\\
&\leq&\frac{3M}{2}\int_{D}|\nabla(u-v)|^2dx+\frac{M}{2}\int_{D}(u-v)^2dx.
\end{eqnarray*}
 Then we have
\begin{eqnarray*}
\frac{1}{2}\frac{d}{dt}\int_{D}|\nabla(u-v)|^2dx &\leq&
C_3\int_{D}|\nabla(u-v)|^2dx+C_4\int_{D}(u-v)^2dx,
\end{eqnarray*}
where $C_3$ and $C_4$ are the constants depending on $C$ and $M$.
By the Gronwall inequality \cite{Evans}, we have
$$
||\nabla(u-v)||_{H^1}^2\leq ||\nabla(g_u-g_v)||_{H^1}^2 \exp(C_2t).
$$
$\Box$

Similar to \cite{MC2}, we have

\textbf{Proof of corollary \ref{cor2}.}

\noindent Since
$$
\frac{1}{2}\frac{d}{dt}\int_{D}|\nabla
u|^2dx=-\int_{D}(u_t)^2dx-\frac{1}{p+1}\frac{d}{dt}\int_{D}u^{p+1}dx
+\frac{1}{2}\frac{d}{dt}\int_{D}au^{p+1}dx,
$$
we have
\begin{equation}\label{energe2}
\lambda(t)+2\int^t_0\int_{D}|u_t|^2dxdt=\int_{D}(|\nabla
g|^2+\frac{2}{p+1}ag^{p+1}-ag^2)dx
+\frac{p-1}{p+1}\int_{D}au^{p+1}dx.
\end{equation}
By the arguments in theorem \ref{thm2}, we have $\lambda(t)$ is
continuous, uniformly bounded in $t\in [0,\infty)$. Moreover, $u\in
L^{\infty}(\mathbb{R}_{+},H^1(D))$ and $u\in
L^{\infty}(\mathbb{R}_{+},L^{p+1}(D))$. Then we can take a
subsequence $\{t_i\}$ with $t_i\to \infty$ such that
$u_i(x)=u(x,t_i)$, $\lambda(t_i)\to\lambda_{\infty}$. By
(\ref{energe2}) and theorem \ref{thm2}, we have
\begin{equation*}
\left\{
\begin{array}{l}
         u_i\to u_{\infty}\ \ \ \ \text{in}\ L^2(D), \\
         u_i\rightharpoonup u_{\infty} \ \ \ \text{in}\ H^1(D)\ \text{and}\ L^p(D),\\
         \partial_t u_i-(\lambda(t_i)-\lambda_{\infty})u_i\to 0 \ \ \ \text{in}\ L^2(D).\\
\end{array}
\right.
\end{equation*}
Since $\partial_t u_i-(\lambda(t_i)-\lambda_{\infty})u_i=\Delta
u_i+\lambda_{\infty}u_i+a(u-u_i^p)$, $u_i\in H^1$ solves the
equation $\Delta
u_{\infty}+\lambda_{\infty}u_{\infty}+a(u_{\infty}-u_{\infty}^p)=0$
in $M$ and $\int_{D}|u_{\infty}|^2dx=1$.
 $\Box$

\section{Conclusion}
In general, we may assume that $a$ is a smooth function both in
space variable and time variable and $a=a(x,t)$ is a periodic
function in time variable $t$. We may also assume the spatial
domain $D$ is a compact manifold (with or without boundary) as in
the works \cite{A98} and \cite{LY}. We leave this subject for
future research.

Based on our study of (\ref{eq2}) above we would like to point out
that the non-local population model of logistic type has the
advantage that the parameter $\lambda(t)$ plays a role like a
control term so that the flow exists globally and has nice
behavior at time infinity. This research shows that global terms
in population modeling should be considered in the future.


\begin{thebibliography}{20}
\bibitem{A98}
T. Aubin, \emph{Some Nonlinear Problems in Riemannian Geometry},
Springer Monogr. Math., Springer-Verlag, Berlin, 1998.

\bibitem{CL09}
C.Caffarelli, F.Lin, \emph{Nonlocal heat flows preserving the
$L^2$ energy}, Discrete and continuous dynamical systems. 23,
49-64 (2009).



\bibitem{DM} Xianzhe Dai and Li Ma, \emph{Mass under Ricci flow},
Commun. Math. Phys., 274, 65-80 (2007).

\bibitem{MD} Y.Du and L.Ma, \emph {Logistic type equations on $\mathbb{R}^{N}$
by a squeezing method involving boundary blow-up solutions},
J.London Math. Soc., 64(2001), 107-124.


\bibitem{Evans} L.Evans,
\emph{Partial Differential Equations,} Graduate studies in Math.,
AMS, 1986

\bibitem{Jiang} Jiang Jifa, Shi, Junping,
\emph{Bistability dynamics in some structured  ecological models},
a chapter in Spatial Ecology, a book published by CRC Press, 2009

 \bibitem{K95}
Kingsland, Sharon (1995). \emph{Modeling Nature: Episodes in the
History of Population Ecology}. University of Chicago Press. pp.
127-146.

\bibitem{MA} L.Ma, A.Q.Zhu, \emph {On a length preserving curve
flow.} Preprint,2008.

\bibitem{MC} L.Ma, L.Cheng, \emph {A non-local area preserving curve
flow.} Preprint, 2008.

\bibitem{MC2} L.Ma, L.Cheng, \emph {non-local heat flows and gradient estimates on closed manifolds.}
J. Evol. Eqs., 2009.


\bibitem{LM}
Li Ma, \emph{Gradient estimates for a simple elliptic equation on
complete non-compact Riemannian manifolds}, Journal of Functional
Analysis, 241(2006)374-382.

\bibitem{LY}
P. Li, S.T. Yau, \emph {On the parabolic kernel of the
Schr\"{o}inger operator}, Acta Math. 156 (1986) 153-01.


\bibitem{R91}
Renshaw, Eric (1991). \emph{Modeling Biological Populations in
Space and Time}. Cambridge University Press. pp. 6-9.


\bibitem{shi09}
Shi, Junping, Sze-Bi Hsu, \emph{Relaxation oscillator profile of
limit cycle in predator-prey system}.
       Discrete and Continuous Dynamical Systems, Series  B, 11 (2009) no. 4, 893-911.


\bibitem{U04}
 Uyenoyama, Marcy; Rama Singh, Ed. (2004).
 \emph{The Evolution of Population Biology}. Cambridge University Press. pp. 1-19.

\end{thebibliography}
\end{document}